\documentclass[10pt,twoside]{article}
\usepackage{graphicx}
\usepackage{amsmath}

\newcommand{\PSbox}[1]{\includegraphics[0in,0in][2.5in,2.5in]{#1}}

\usepackage{amssymb}

\usepackage{Latex-document}

\title{\bf Quantum Resonances  and \vskip -2mm
Partial Differential Equations \vskip 6mm}

\author{M. Zworski\vspace*{-0.5cm}\thanks{Department of Mathematics,
University of California, Berkeley, CA 94720, USA. E-mail:
zworski@math.berkeley.edu}}
\date{\vspace{-8mm}}

\begin{document}
\maketitle

\thispagestyle{first} \setcounter{page}{243}

\begin{abstract}\vskip 3mm

Resonances, or scattering poles, are complex numbers which
mathematically describe meta-stable states: the real part of
a resonance gives the rest energy, and its imaginary
part, the rate of decay of a meta-stable state.
This description emphasizes the quantum mechanical aspects of
this concept but similar models appear in many
branches of physics, chemistry and mathematics, from molecular dynamics
to automorphic forms.

In this article we will will describe the recent progress  in the
study of resonances based on the theory of partial differential equations.

\vskip 4.5mm


\end{abstract}

\vskip 12mm

\section{Introduction} \label{section 1}\setzero
\vskip-5mm \hspace{5mm }

Eigenvalues of self-adjoint operators appear naturally in quantum
mechanics and are in fact what we observe experimentally in many
situations. To explain the need for the more subtle notion of
{\em quantum resonances} we consider the following simple example.

\begin{figure}[htb]
\begin{center}
\label{f.1}
 \vspace{0.3cm}
   %\hspace{0.2cm}
   \includegraphics{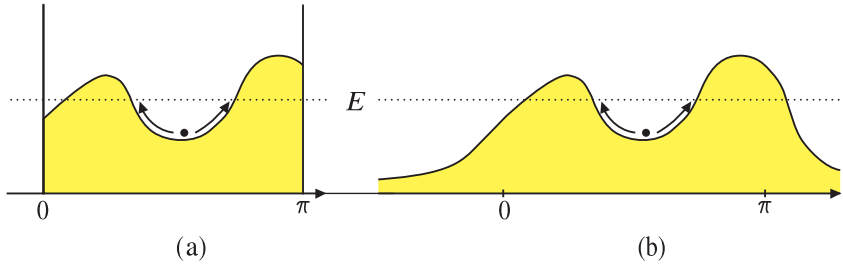}
%   \includegraphics{rpg.eps}
%   \vspace{-0.8cm}
%\Dessin{rpg2}
\begin{minipage}[h]{10cm}
\caption{{(a)} A potential well on a finite interval. {(b)} A
potential on the whole line the same classical picture but a very
different quantum picture. }
\end{minipage}
\end{center}
\end{figure}

Let $ V ( x ) $ be  a potential  on a bounded interval,
as shown in Fig.\ref{f.1} a). If $ \xi $ denotes the classical
momentum, then the classical energy is given by
\begin{equation}
\label{eq:cham}
E = \xi^2 + V ( x) \,,
\end{equation}
and the motion of a classical particle can deduced from this
equation by considering $ E $ as the Hamiltonian. The quantized
Hamiltonian is given by
\begin{equation}
\label{eq:qham}
P( h ) = ( h D_x )^2 + V ( x) \,, \ \ \xi \mapsto h D_x \,, \ \ D_x=
\frac{1}{i} \partial_x \,.
\end{equation}
The operator $ P ( h ) $ considered as an unbounded operator on $ L^2 $
of a bounded interval (with, say, Dirichlet boundary conditions)
has a discrete spectrum, which gets denser and
denser as $ h \rightarrow 0 $, which corresponds to getting closer to
classical mechanics:
\begin{equation}
\label{eq:heig}
P ( h ) u ( h ) =  E ( h ) u ( h ) \,, \ \ E ( h ) \in {\bf R} \,, \ \
\int |u ( h ) |^2 < \infty \,, \ \ u ( h ) ( \pm \pi ) = 0 \,.
\end{equation}
 The eigenvalues, $ E( h) $, are what we (in principle) observe,
and the square-integrable eigenfunctions, $ u ( h )$,  are the wave functions.

Now, consider the same example but with a potential $ V( x )$ on
$ {\bf R} $, as shown in Fig.\ref{f.1} b). At the energy $ E $, and inside
the well created by the potential, the classical motion is the
same as in the previous case. Hence,
we expect that (at least for $ h $ very small) there should exist
a quantum state corresponding to the classical one. It is clear in
dimension one that for the potential shown in the figure, $ P ( h ) $
given by \eqref{eq:qham} has {\em no} square integrable eigenfunctions
\eqref{eq:heig}.

\vskip -3mm

A non-obvious remedy for this is to think of eigevalues (in the
case of Fig.\ref{f.1}a) as the poles of the resolvent of $ P ( h )
$:
\[ R ( z, h ) = ( P ( h ) - z )^{-1} \,.\]
In the case of Fig\ref{f.1}b), the resolvent is not bounded on $ L^2 (
{\bf R} ) $ for $ z > 0 $, which corresponds to free motion or
scattering. However,
under suitable assumptions on $ V ( x ) $ near infinity we have
a meromorphic continuation of $ R ( z, h ) $ from $ {\rm{Im}} z > 0 $ to the
lower half-plane:
\[ R ( z , h ) \; : \; L^2_{\rm{comp}} ( {\bf R} ) \longrightarrow
 L^2_{\rm{loc}} ( {\bf R} )\,. \]
The complex poles of $ R ( z , h ) $ are the replacement for eigenvalues,
are called {\em resonances}. For a recent presentation of
these in the physical
literature and for examples of current interest see for instance \cite{G2}.

In this review, we will restrict ourselves to the case of
\[  P ( h ) = -h^2 \Delta + V ( x) \,, \ \ x \in {\bf R}^n \,, \]
where the potential $ V ( x) $ satisfies some decay and analyticity
assumptions near infinity, or is simply compactly supported.
Then our discussion above applies and resonances are defined as
the poles of $ R ( z, h ) $ -- see \cite{lm} for an attractive
new presentation of the meromorphic continuation properties and
references. In our convention the resonances are located in the lower
half-plane.

We should stress that most of the results hold for
very general {\em black box perturbations} of Sj\"ostrand-Zworski \cite{sz}
which allow for a study of diverse problems without going into their
specific natures -- see \cite{s} for definitions.
Also, despite the fact that the motivation presented here
came from molecular dynamics,
similar issues arise in other settings, from
automorpic scattering to
electromagnetic (obstacle)
scattering  -- see \cite{GuZw},  and \cite{sz7} respectively.

We suggest the the surveys \cite{Mel-b},\cite{s},\cite{Vo},\cite{Zws}, and
\cite{ZwS} for the review of earlier results, and we concentrate
on the progress achieved in the last few years. For an account of
work based on other methods and motivated by different physical
phenomena we refer to \cite{MeSi}.

\section{Trace formul{\ae} for resonances} \label{section 2}
\setzero\vskip-5mm \hspace{5mm }

Trace formul{\ae} provide one of the most elegant descriptions of the
classical-quantum correspondence. One side of a formula is given
by a trace of a quantum object, typically derived from a quantum
Hamiltonian, and the other side is described in terms of closed
orbits of the corresponding classical Hamiltonian.
 A new general approach based on {\em quantum
monodromy operator} which quantizes the Poincar\'e map in a natural
way was recently given in \cite{SZ8}.

In general the spectral (or scattering) side of the formula is given
in terms of the trace of a function of the quantum Hamiltonian
$ f ( P ( h ) ) $. In the case of self-adjoint problems with discrete
sets of eigenvalues  the spectral theorem readily provides an
expression for $ {\rm{tr}} \; f (  P ( h )) $ but the problem becomes
subtle for resonances. It has been studied by Lax-Phillips,
Bardos-Guillot-Ralston, Melrose, Sj\"ostrand-Zworski, Guillop\'e-Zworski
-- see \cite{z1} and references given there.

More recently Sj\"ostrand \cite{s},\cite{Sj-trace} introduced
{\em local trace formul{\ae}} for resonances, and that concept was
developed further in \cite{pz3} and \cite{pb}. For $ P ( h ) $ given
by \eqref{eq:qham} with $ V ( x) $ decaying sufficiently fast at
infinity it can be stated as follows.

Let $ \Omega $ be an open, simply connected, pre-compact
subset of $ \{ {\rm{Re}}\; z > 0 \} \subset {\bf C} $, such that
$ \Omega \cap {\bf R} $ is connected. Suppose
that $ f $ is holomorphic in a neighbourhood of $ \Omega $ and
that $ \psi \in C^\infty_{\rm{c}} ( {\bf R}) $ is equal to $ 1 $ in
$ \Omega \cap {\bf R} $ and supported in a neighbourhood of $ \Omega \cap
{\bf R} $. Then, denoting the resonance set of $ P ( h ) $ by $ {\rm{Res}}
( P ( h ) )$,
\begin{gather}
\label{eq:trace}
\begin{gathered}
{\rm{tr}}\; \left( ( f \psi ) ( P ( h ) ) - ( f \psi ) ( -h^2 \Delta )
\right) = \sum_{ z \in {\rm{Res}} ( P ( h )) \cap \Omega } f( z )
+ E_{f,\psi}( h)  \,,\\
| E_{f , \psi } ( h )
|\leq h^{-n} C_{ \Omega, \psi}  \max \{ |f ( z ) |\; : \;
z \in \overline{ \Omega}_1 \setminus \Omega\,, \ {\rm{Im}} \;  z \leq 0  \} \,,
\end{gathered}
\end{gather}
where $ \Omega_1 $ is a neighbourhood of $ \Omega $ (see \cite{pb}
and \cite{pz3} for more
precise versions).

The basic upper bound on the number of resonances is given as follows:
\begin{equation}
\label{eq:upper}
\sharp \; {\rm{Res}}(P(h)) \cap \Omega = {\mathcal O} ( h^{-n}) \,,
\end{equation}
see \cite{Mel-b},\cite{Sj-d},\cite{s},\cite{Zws}. It would consequently
appear that the error term in \eqref{eq:trace} is of the same order as the
sum over the resonances. However,
by choosing the function $ f $ so that it is small in $  \Omega_1
\setminus \Omega \cap \{ {\rm{Im}}\; z \leq 0 \} $, the sum of $ f ( z ) $
can dominate the left hand side of \eqref{eq:trace}. Doing that,
Sj\"ostrand has shown that an analytic singularity of
$ E \mapsto | \{ x \; : \; V ( x) \geq E \} | $ at $ E_0 $ gives a lower bound
in \eqref{eq:upper} for $ \Omega $ a neighbourhood of $ E_0 $ -- see
\cite{Sj-trace} and references given there.

Another application of local trace formul{\ae} techniques is in analysing
resonances for {\em bottles}, that is perturbations of the Euclidean space
in which the ``size'' of the perturbation may grow but which are connected
to the euclidean infinity through a fixed ``neck of the bottle'' -- see
\cite{Sj-trace},\cite{pz3}.

\section{Breit-Wigner approximations} \label{section bw}
\setzero\vskip-5mm \hspace{5mm }

In a scattering experiment the physical data is mathematically encoded
in the scattering matrix, $ S ( \lambda , h ) $. It is in the behaviour of
objects derived from the scattering matrix that we see physical manifestations
of abstractly defined resonances. The classical and still central
Breit-Wigner approximation provides this connection. The most mathematically
tractable case is provided by the {\em scattering phase} $ \sigma
( \lambda , h )
= \log \det S ( \lambda, h  ) / ( 2 \pi i ) $, which is defined for $ V $'s
with sufficient decay (otherwise relative phase shifts have to be used,
see \cite{pb}). When $ \lambda $ is close to a resonance $ E - i \Gamma $,
the Breit-Wigner approximation says that
\begin{equation}
\label{eq:bw1}
\sigma' ( \lambda , h ) \sim \frac{1}{\pi} \frac{  \Gamma }{ ( \lambda - E)^2
+ \Gamma^2 } \,,
\end{equation}
that is, if $ \Gamma $ is small,
$ \sigma ( \lambda ) $ should change by approximately $ 1 $ as
$ \lambda $ crosses $ E $. This has been justified rigorously in
some situations in whch a given resonance, close to the real axis, is isolated.

In view of \eqref{eq:upper} the number of resonances in fixed regions
can be very large and for $ h$ small clouds of resonances need to be
considered to obtain a correct form of the Breit-Wigner approximation.
A formalism for that was introduced by Petkov-Zworski \cite{PZ},\cite{pz3},
and it was developed further by Bruneau-Petkov \cite{pb} and
Bony-Sj\"ostrand \cite{BS}. It is closely related to extending the
trace formula \eqref{eq:trace} to $h$ dependent $ \Omega $'s.
To formulate it, let us introduce $ \omega ( z , E ) = (1/\pi)
\int_E ( |{\rm{Im}}\;
 z | / | z- \lambda|^2 ) d\lambda $, $ {\rm{Im}}\; z  <  0 $,
the harmonic measure corresponding to the upper half-plane.
Then, for non-critical $ \lambda$'s (that is for $ \lambda$'s for which
$ \xi^2 + V ( x) = \lambda  \Rightarrow d_{( x , \xi)  } ( \xi^2 + V( x) )
\neq 0 $)
\begin{equation}
\label{eq:bw2}
\sigma ( \lambda + \delta , h ) -
\sigma ( \lambda - \delta , h )  =
\sum_{{ z \in {\rm{Res}}( P ( h )) } \atop { | z - \lambda | <  h } }
\omega ( z , [ \lambda - \delta , \lambda + \delta ]) +
{\mathcal O} ( \delta ) h^{-n} \,,
\end{equation}
$  0 < \delta < h /C $. The main point is that $ \delta $ can be made
arbitrarily small and that we only need to include resonances close to
$ \lambda $. If we do not assume that $ \lambda $ is non-critical
weaker, but still useful, results can be obtained from factorization of
the scattering determinant \cite{pz3}, or more generally, from the
the analysis of the phase shift function \cite{pb}.

One of the consequencies of the development of the Breit-Wigner
approximations for clouds of resonances were new estimates on the
number of resonances in small regions, first in \cite{PZ}, and then
in greater generality in \cite{B},\cite{pb},\cite{pz3}.
If in place of \eqref{eq:upper} we had a Weyl law with a remainder
$ {\mathcal O} ( h^{1-n} ) $, then, as for eigenvalues,
\begin{equation}
\label{eq:supper}
\sharp \; \{ z \in {\rm{Res}} (P( h ) ) \; : \; | z - \lambda | <  \delta \}
= {\mathcal O} ( \delta ) h^{-n} \,, \ \  C h < \delta < 1/C \,.
\end{equation}
It turns out that despite the lack of the Weyl law we still have
this estimate for non-critical $ \lambda $'s.

\section{Resonance expansions of propagators} \label{section re}
\setzero\vskip-5mm \hspace{5mm }

In the case of discrete spectrum of a self-adjoint
operator the propagator, $ \exp ( -i t P ( h ) /h ) $ can be expanded
in terms of the eigenvalues. In fact, our understanding of
``state specificity'' often comes from such ``Fourier decompositions'' into
modes. For problems in which decay or escape to infinity are possible
such expansions are still expected but just as in the case of trace
formul{\ae} far from obvious. For non-trapping perturbations and in
the context of the wave equation the expansions in terms of
resonances were studied in late 60's by Lax-Phillips and Vainberg
(see the references in \cite{tz2}).

For trapping perturbations the expansions were investigated by
Tang-Zworski \cite{tz2}, and then by Burq-Zworski \cite{bz},
Christiansen-Zworski \cite{cz}, Stefanov \cite{St},
and most recently by Nakamura-Stefanov-Zworski \cite{NSZ}.

Here, we will recall the general expansion given in \cite{bz}: let
$ \chi \in C^\infty_{\rm{c}} ( {\bf R}^n ) $, $ \psi \in
C^\infty_{\rm{c} }( 0 , \infty) $, and let  $ {\rm{chsupp}}\; \psi
= [a, b] $. There exists $  0 < \delta < c ( h ) < 2 \delta $ and
$ L $, so that we have
\begin{equation}
\label{eq:bz}
\begin{split}
\chi e^{ - i t  P (h)  / h } \chi \psi ( P ( h )  ) = &
\sum_{ {z \in \Omega (h) \cap \rm{Res}(P)} }
 \chi {\rm{Res}} ( e^{-  i  t \bullet /h } R ( \bullet, h ), z ) \chi
\psi ( P ( h ) ) \\
\; &  + \; {\mathcal O}_{ L^2  \rightarrow L^2 } ( h^\infty ) \,, \ \ \
{\text{for}} \ \
t > h^{-L }\,,
\end{split}
\end{equation}
$\Omega ( h ) = ( a - c(h), b + c(h) ) - i [ 0 , 1/C ) \,,$
and where $ {\rm{Res}} ( f ( \bullet ) , z ) $ denotes the residue of
a meromorphic family of operators, $ f $, at $ z $.

The function $ c(h) $ depends on the distribution of
resonances: roughly speaking we cannot ``cut'' through a dense
cloud of resonances. Even in the very well understood case of
the modular surface \cite[Theorem 1]{cz} there is, currently at
least, a need for some non-explicit grouping of terms.
This is eliminated by the separation condition \cite[(4.4)]{tz2}
which however is hard to verify.

The unpleasant feature of \eqref{eq:bz} is the need for very
large times $ \sim h^{-L }$ and the presence of
a non-universal parameter, $ c ( h ) $.
The former is necessary in this formulation as one sees
by considering the free case $ P ( h ) = - h^2 \Delta $.
However an expansion valid for all times is possible
in the case when a part of $ V $ constitutes a {\em barrier}
separating the trapped set in  $ \{ ( x, \xi ) \; : \;
\xi^2 + V ( x ) \in  {\rm{supp}} \; \psi \} $ from infinity \cite{NSZ}.
The example shown in Fig.\ref{f.1} b) is of that type.
By the trapped set in $ \Sigma \subset T^* {\bf R}^n $ we mean the set
\begin{equation}
\label{eq:trapped}
 K \cap \Sigma = \{ ( x , \xi ) \in \Sigma \; : \;
| \exp ( t H_p ) ( x , \xi ) | \not \longrightarrow \infty \,, \
t \longrightarrow \pm \infty \} \,, \end{equation}
where $ H_p $ is the Hamilton vectorfield of $ p = \xi^2 + V ( x ) $.
One of the  components comes from the work of Stefanov \cite{St-ajm}
on making estimates \eqref{eq:upper} more quantitative with constants
related to the volume of the trapped set (for yet finer results of
that type see Section.6).

\section{Separation from the real axis} \label{section sr}
\setzero\vskip-5mm \hspace{5mm }

One of the most striking applications of PDE techniques in the
study of resonances is the
work of Burq (see \cite{b2} and references given there) on
the separation of resonances from the real axis, estimates of
the cut-off resolvent on the real axis, and the consequent
estimates on the time decay of energy.

For non-trapping perturbations we have the following estimate
on the truncated meromorphically continued resolvent
\begin{equation}
\label{eq:res1}
\| \chi R ( z , h ) \chi \|_{ L^2 \rightarrow L^2  } \leq
C \frac{ \exp({ C|{\rm{Im} }\; z | / h }) }{h} \,, \ \ {\rm{Im}} \; z
> - M h \log \frac{1}{h} \,, \ \ \chi \in C^\infty_{\rm{c}} ( {\bf
R}^n ) \,,
\end{equation}
for any $ M$, see \cite{M} for the absence of resonances (implicit in
\eqref{eq:res1}), and \cite{NSZ} for the (easily derived) estimate.

Since the work of Stefanov-Vodev and Tang-Zworski (see
\cite{St-ajm},\cite{tz},\cite{Vo} and references given there) we know that
in many trapping situations\footnote{roughly speaking, whenever there exists
an elliptic orbit of the $ H_p$-flow}, there exist many resonances
converging to the real axis\footnote{That in specific trapping
situations there exist resonances whose distance to the real axis is
of order $ \exp( - S/ h ) $ is classical in physics with the most
precise mathematical results given by Helffer-Sj\"ostrand -- see
\cite{lm}. Here we concentrate on very general existence results which
guarantee lower bounds on the number of resonances.}.
The question then is how close could the resonances approach the real
axis or, as it turns out equivalently, how big can the truncated resolvent
be on the real axis. Heuristically, tunneling should prevent arbitrary
closeness to the real axis, and the separation should be universally
given by at least $ \exp ( -S/ h ) $.

Tunneling of solutions is made quantitative in PDEs through {\em Carleman
estimates} used classically to show unique continuation of solutions of
second order equations. Following the work of Robbiano, and Lebeau-Robbiano,
Burq succeeded in applying Carleman estimates to a wide range of
resonance problems -- see \cite{b2} and references given there.
In the setting described here his basic result says that if $
\Omega $ is as in \eqref{eq:trace} then
\begin{gather}
\label{eq:burq}
\begin{gathered}
\exists \; S_1 , S_2 \,, \ \
z \in {\rm{Res}} ( P ( h ) ) \cap \Omega \; \Longrightarrow \;
{\rm{Im}} \;  z > - \exp ( - S_1 / h ) \,,\\
\| \chi R ( z, h ) \chi \|_{ L^2 \rightarrow L^2  } \leq
\exp ( S_2 / h ) \,, \ \ z \in \Omega \cap {\bf R}\,,
\ \ \chi \in C^\infty_{\rm{c}} ( {\bf R}^n ) \,.
\end{gathered}
\end{gather}
In addition, improved estimates are possible if $ \chi $ is assumed
to have support outside the projection of the trapped set:
\begin{gather}
\label{eq:res2}
\| \chi R ( z, h ) \chi \|_{ L^2 \rightarrow L^2  } \leq
C/ h  \,, \ \ z \in \Omega \cap {\bf R}\,,
\ \ \chi \in C^\infty_{\rm{c}} ( {\bf R}^n \setminus \pi ( K )  ) \,.
\end{gather}
That means that away from the interaction region $ \pi ( K ) $ we
have, on the real axis, the same estimates as in the non-trapping case
\eqref{eq:res1} and that has immediate applications in scattering theory
\cite{St-sc}.

\section{Resonances in chaotic scaterring}
\setzero\vskip-5mm \hspace{5mm }

Since the work of Sj\"ostrand \cite{Sj-d} on geometric upper bounds for the
number of resonances, it has been expected that for chaotic scattering
systems the density of resonances near the real axis can be approximately
given by a power law with the power equal to half of the dimension of the
trapped set (see \eqref{eq:conj} below). Upper bounds in geometric
situations have been obtained in \cite{WZ} and \cite{ZwI}.

\begin{figure}
\label{f:2}
  \begin{center}
  \PSbox{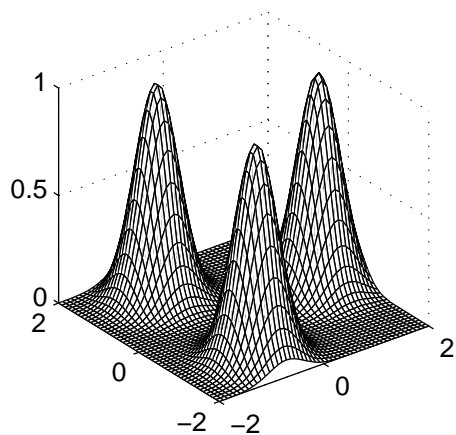}
%  \PSbox{pot.eps}
\begin{minipage}[h]{10cm}
  \caption{Graph of a potential $ V ( x) $, for which the classical flow
of the Hamiltonian $ \xi^2 + V ( x) $ is hyperbolic for energies close to
$ 0.5$.}
\end{minipage}
\end{center}
\end{figure}

An example of a potential in $ V \in C^\infty_{\rm{c}} ( {\bf R}^2
) $ for which the flow of $ H_p $, $ p = \xi^2 + V ( x ) $, is
hyperbolic (and hence scattering exhibits chaotic features) is
shown in Fig.\ref{f:2}

\begin{figure}
\label{f:3}
  \begin{center}
  \PSbox{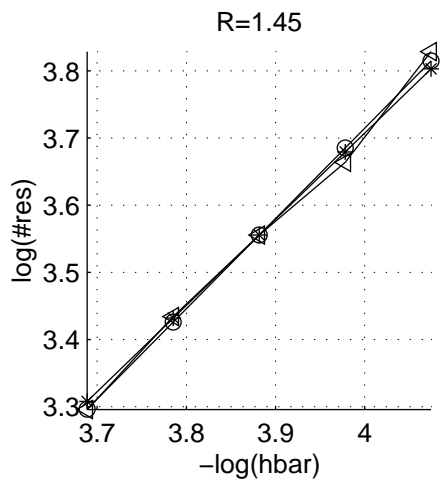}
%  \PSbox{hist1.eps}
\begin{minipage}[h]{10cm}
  \caption{A plot of $\log(N_{res})$ as a function for $-\log(\hbar)$
  for the potential shown in in Fig. \ref{f:2}
  The value of $\hbar$ ranges from $0.025$ to $0.017$.
  Triangles represent numerical data, circles least squares regression,
  and stars the slope predicted by the conjecture.}
\end{minipage}
  \end{center}
\end{figure}

A recent numerical study \cite{LZ} for that potential
indicates that the density of resonances satisfies a lower
bound related to the dimension of the trapped set -- see Fig.\ref{f:3}
In a complicated semi-classical situation
studied in \cite{LZ}, the dimension is a delicate
concept and it may be that different notions of dimension have to be
used for upper and lower bounds. That point is emphasized in \cite{lsz}
where numerical data for semi-classical zeta function for several
convex obstacles is analyzed.

In the case of {\em convex co-compact hyperbolic quotients}, $ X =
\Gamma \backslash {\bf H}^2 $, studied in \cite{ZwI} the situation
is particularly simple as the quantum resonances coincide with the
zeros of the zeta-function -- see \cite{PP}. The notion of the
dimension of the trapped set is also clear as it is given by $ 2 (
1 + \delta ) $. Here $ \delta = \dim \Lambda( \Gamma ) $ is the
dimension of the limit set of $ \Gamma$, that is the set of
accumulation points of the elements of $ \Gamma $ (they are all
hyperbolic), $ \Lambda ( \Gamma ) \subset \partial {\bf H}^2 $.

Hence we expect that
\begin{equation}
\label{eq:conj}
  \sum_{ |{\rm{Im}}\;  s | \leq r \,, \; {\rm{Re}} \; s > - C } m_\Gamma
 ( s ) \sim r^{ 1+ \delta}
\, , \end{equation}
where $ m_\Gamma
 ( s ) $ is the multiplicity of the zero of the zeta function of
$ \Gamma $ at $ s $. We also changed the traditional convention: here $
h^2 s ( 1 -s ) = z $ in the notation of previous sections, and
$ h (  {\rm{Re}}\; s - 1/2 ) $, $ h  {\rm{Im}}\; s  $ correspond to
$ {\rm{Im}}\; z  $, $  {\rm{Re}}\; z $ respectively.

An upper bound of this form was established in \cite{ZwI} but
a simpler method giving improved upper bounds has been recently
presented in \cite{zrs}. The new method is based on zeta function
techniques. What we obtain is a bound in the case of convex co-compact
Schottky groups
\begin{equation}
\label{eq:ubd}
  \sum \{
m_\Gamma (s) \; : \;  {r \leq {\rm{Im}}\;  s  \leq r + 1}\,, \ \
 { {\rm{Re}}\;  s  > - C_0 }
\}  \leq  C_1 r^{ \delta}
\, , \end{equation}
where $ \delta = {\rm{dim}} \; \Lambda ( \Gamma ) $.
This improved estimate is a ``fractal'' version of \eqref{eq:supper}:
$ 1 + \delta $ now plays the r\^ole of $ n $. It is clear that the
method works in greater generality and implementing it is part of an
ongoing project.

\label{lastpage}

\end{document}